
\input  math.macros
\input  Ref.macros
\input EPSfig.macros
\input labelfig.tex

\checkdefinedreferencetrue
\theoremcountingtrue
\sectionnumberstrue
\figuresectionnumberstrue
\forwardreferencetrue
\citationgenerationtrue
\hyperstrue
\nobracketcittrue
\initialeqmacro

\bibsty{myapalike}
\bibfile{\jobname}

\long\def\comment#1{\relax}

\def\subsection#1{\bigskip\noindent {\bf #1.}\par\nobreak}

\def\I#1{{\bf 1}_{#1}}
\def\II#1{\I{\{#1\}}}

\let\circ=\diamond
\let\hat=\widehat

\def\Zd{\Z^d}
\def\vi{v_R^*}
\def\Rel#1{{\cal #1}}
\def\RR{\Rel R}
\def\PP{\Rel P}
\def\LL{\Rel L}

\def\UU{\Rel U}
\def\MM{\Rel Q}

\def\Y{\Upsilon}

\def\field#1{{\frak #1}}
\def\dd{{\rho}}
\def\st{:\,}

\def\o#1{\langle#1\rangle}
\def\oo#1#2{\o{#1}^{\,#2}}
\def\ta#1/{\o{#1}}
\def\sm{\smallsetminus}
\def\lec{\preccurlyeq}
\def\gec{\succcurlyeq}
\def\C{{\cal C}}
\def\mmS{{\cal S}}
\def\Q{Q}
\def\v{v_0}
\def\A{A}
\def\aB{\overrightarrow{B}}

\def\xx{x}
\def\gar{\gamma_R}
\def\AR{A_R}
\def\mid{|}

\vglue20pt

\title{Geometry of the Uniform Spanning Forest:}
\title{Transitions in Dimensions 4, 8, 12, \dots }




\author{{Itai Benjamini, Harry Kesten, Yuval Peres, Oded Schramm}}

\abstract{%
The uniform spanning forest (USF) in $\Zd$
is the weak limit of random, uniformly chosen,
spanning trees in $[-n,n]^d$.
\ref B.Pemantle/ proved that the USF consists a.s.\
of a single tree if and only if $d \le 4$.
We prove that any two components of the USF in $\Z^d$ are adjacent
a.s.\ if $5 \le d \le 8$, but not if $d \ge 9$.
More generally,
let $N(x,y)$ be the minimum number of edges outside the USF
in a path joining $x$ and $y$ in $\Zd$.
Then
$$
\max\bigl\{N(x,y) \st x,y\in\Zd\bigr\} = \bigl\lfloor (d-1)/4
\bigr\rfloor
\hbox{ a.s. }
$$
The notion of stochastic dimension for random relations
in the lattice is introduced
and used in the proof.
}

\bottom{Primary
60K35
. %
Secondary
60J15
.}
{Stochastic dimension
, Uniform spanning forest.}
{Research partially supported
by NSF grants DMS-9625458 (Kesten) and DMS-9803597 (Peres), and
by a Schonbrunn Visiting Professorship (Kesten).}

\vfill\eject

\bsection {Introduction}{s.intro}

A {\bf uniform spanning tree} (UST) in a finite graph
is a subgraph chosen uniformly at random among all
spanning trees.  (A spanning tree is a subgraph
such that every pair of vertices in the original graph
are joined by a unique simple path in the subgraph.)
The {\bf uniform spanning forest} (USF) in $\Zd$ 
is a random subgraph of $\Zd$, that was defined by
\ref b.Pemantle/ (following a suggestion of R. Lyons), as follows:
The USF is the weak limit of  uniform spanning trees in larger and
larger finite boxes.   
Pemantle showed that the limit exists, that it
does not depend on the sequence of boxes, and
that every connected 
component of the USF is an infinite tree.
 See Benjamini, Lyons, Peres and Schramm (2001)
(denoted \ref b.BLPS:usf/ below) for a thorough
study of the construction and properties of the
USF, as well as references to other works on the subject.
Let $T(x)$ denote the tree in the USF which contains the vertex
$x$. Also define
$$
N(x,y)
= \min\bigl\{\hbox{number of edges outside the USF in a path from }
              x \hbox{ to } y\bigr\}
$$
(the minimum here is over all paths in $\Zd$ from $x$ to $y$).

\ref b.Pemantle/ proved that
for $d \le 4$,  almost surely $T(x)= T(y)$ for all $x,y\in \Zd$, and
for $d>4$, almost surely $\max_{x,y} N(x,y)>0$.
The following theorem shows that a.s.\ $\max N(x,y)\le 1$ for
$d=5,6,7,8$,
and that $\max N(x,y)$ increases by
$1$ whenever the dimension $d$ increases by $4$.

\procl t.main1
$$
\max\bigl\{ N(x,y)\st x,y\in\Zd\bigr\}
= \left\lfloor {{d-1}\over 4} \right\rfloor \qquad\hbox{ a.s.}
\label e.main1
$$
Moreover, a.s.\ on the event $\bigl\{T(x) \ne T(y)\bigr\}$, there exist
infinitely many disjoint simple paths in $\Zd$ which
connect $T(x)$ and $T(y)$ and which contain at most
$\lfloor (d-1)/4 \rfloor$ edges outside the USF.
\endprocl

It is also natural to study
$$
D(x,y):= \lim_{n \to \infty} \inf\{|u-v|: u \in T(x), v \in
T(y), |u|, |v| \ge n\},
$$
where $|u|=\|u\|_1$ is the $\l^1$ norm of $u$.
The following result is a consequence of \ref b.Pemantle/ and our
proof of \ref t.main1/.

\procl t.main2 
Almost surely, for all $x, y \in \Zd$,
$$
D(x,y) =\cases{0 &if  $d\le 4$,\cr
               1 &if $5 \le d \le 8$ and $T(x)\neq T(y)$,\cr
               \infty &if  $d \ge 9$ and $T(x)\neq T(y)$.\cr}
$$
\endprocl

When $5\leq d\leq 8$, this provides a
natural example of a translation invariant
random partition of $\Zd$, into infinitely many
components, each pair of which come infinitely often within unit
distance
from each other.

The lower bounds on $N(x,y)$ follow readily from
standard random walk estimates
(see \ref s.mainproofs/), so the bulk of our work
will be devoted to the upper bounds.

Part of our motivation comes from
the conjecture of \ref B.Newman-Stein/
that invasion percolation clusters in  $\Zd$, $d\ge 6$, are in
some sense 4-dimensional and that two such clusters,
which are formed by starting at two different
vertices,  will intersect with probability 1 if  $d < 8$,
but not if $d>8$.
A similar phenomenon is expected for minimal spanning trees on
the points of a homogeneous Poisson process in $\R^d$.
These problems are still open, as the tools presently available
to analyze invasion percolation and minimal spanning forests
are not as sharp as those available for the uniform
spanning forest.

In the next section the notion of stochastic dimension is introduced.
A random relation  $\RR\subset\Zd\times\Zd$ has
{\bf stochastic dimension} $d-\alpha$, if there is some constant $c>0$
such
that
for all $x\neq z$ in $\Zd$,
$$
c^{-1}|x-z|^{-\alpha}\le\,\P[x\RR z]\le c |x-z|^{-\alpha},
$$
and if a natural correlation inequality \ref e.d/
(an upper bound  for $\P[x\RR z,\, y\RR w]$)
holds.
The results regarding stochastic dimension are formulated and proven
in this generality, to allow for future applications.

\medskip

The bulk of the paper is devoted to the proof of the
upper bound on $\max N(x,y)$ in~\ref e.main1/.
We now present an overview of this proof.
Let $\UU^{(n)}$ be the relation $N(x,y)\le n-1$.
Then $x\UU^{(1)}y$ means that $x$ and $y$ are
in the same USF tree, and $x\UU^{(2)}y$ means that $T(x)$
is equal or adjacent to $T(y)$.  We show that
$\UU^{(1)}$ has stochastic dimension $4$ when $d\ge 4$.
When $\RR,\LL\subset\Zd\times\Zd$ are independent relations with
stochastic dimensions $\dim_S(\RR)$ and $\dim_S(\LL)$, respectively,
it is proven that the composition $\LL\RR$
(defined by $x\LL\RR y$ iff there is a $z$ such that
$x\LL z$ and $z\RR y$) has stochastic dimension
$\min\{\dim_S(\RR)+\dim_S(\LL),d\}$.
It follows that the composition of $m+1$ independent copies
of $\UU^{(1)}$ has stochastic dimension $d$, where 
$m$ is equal to the right hand of~\ref e.main1/.
By proving that $\UU^{(m+1)}$ stochastically dominates the composition
of $m+1$ independent copies of $\UU^{(1)}$, we conclude
that $\dim_S(\UU^{(m+1)})=d$, which implies
$\inf_{x,y\in\Zd}\P[N(x,y)\le m]>0$.
Non-obvious tail-triviality arguments then give
$\P[N(x,y)\le m]=1$ for every $x$ and $y$ in $\Z^d$,
which proves the required upper bound.

\medskip

In \ref s.usfprop/ we present the relevant USF
properties needed; in particular, we obtain a tight upper
bound, \ref t.spread/, for the probability that a finite set of vertices
in $\Zd$ is contained in one USF component.
Fundamental for these results is a method from  \ref b.BLPS:usf/
for generating the USF in any transient graph, which is
 based on an algorithm by
\ref b.Wilson:alg/ for sampling uniformly from the
spanning trees in finite graphs. (We recall this method in \ref
s.usfprop/.)

Our main results are established
in \ref s.mainproofs/. \ref s.examp/ describes several examples
of relations having a stochastic dimension, including long-range
percolation, and suggests some conjectures.
We note that proving $D(x,y)\in\{0,1\}$
for $5 \le d \le 8$, is easier than the higher dimensional result.
(The full power of \ref t.comps/ is not needed; \ref c.low/ suffices).

\bsection {Stochastic dimension and compositions}{s.comp}

\procl d.spread
  When $x,y\in\Zd$,
we write $\o{xy}:=1+|x-y|$, where $|x-y|=\|x-y\|_1$ is
the distance from $x$ to $y$ in the graph metric on $\Zd$.
  Suppose that $W\subset\Zd$ is finite, and $\tau$ is a tree
on the vertex set $W$ ($\tau$ need not be a subgraph
of $\Zd$).  Then let $\ta\tau/:=\prod_{\{x,y\}\in\tau}\o{xy}$
denote the product
of $\o{xy}$ over all undirected edges $\{x,y\}$ in $\tau$.
Define the {\bf spread} of $W$  by
$\ta W/:=\min_{\tau} \ta \tau/$, where
$\tau$ ranges over all trees on the vertex set $W$.
\endprocl

For three vertices,
$
\o{xyz}=\min\{\o{xy}\o{yz},\o{yz}\o{zx},\o{zx}\o{xy}\}
$.
More generally, for $n$ vertices, $\o{x_1 \dots x_n}$
 is a minimum of $n^{n-2}$ products (since this is
the number of trees on $n$ labeled vertices);
see \ref r.spready/ for a simpler
equivalent expression.

\procl d.dimdef \procname{Stochastic dimension}  
Let $\RR$ be a random subset of $\Zd\times\Zd$.
We think of $\RR$ as a relation, and usually write
$x\RR y$ instead of $(x,y)\in\RR$.
Let $\alpha \in [0,d)$.
We say that $\RR$ has {\bf stochastic dimension} $d-\alpha$,
and  write $\dim_S(\RR)=d-\alpha$, if
there is a constant $C=C(\RR)<\infty$ such that
$$
C\, \P[x \RR z] \ge  \o{xz}^{\, -\alpha}
\,,
\label e.a
$$
and
$$
\P[x\RR z,\ y\RR w] \le C\, \oo{xz}{-\alpha}\oo{yw}{-\alpha}+
C\,\ta xzyw/^{-\alpha}
\,,
\label e.d
$$
hold for all $x,y,z,w\in\Zd$.
\endprocl

Observe that \ref e.d/ implies
$$
\P[x\RR z] \le  2C \, \o{xz}^{\, -\alpha},
\label e.dd
$$
since we may take $x=y$ and $z=w$.
Also, note that $\dim_S(\RR)=d$  iff $\inf_{x,z \in\Zd}\P[x \RR z]>0$.

To motivate \ref e.d/, focus on the special case in which $\RR$ is
a random equivalence relation. Then heuristically, the first summand
in \ref e.d/ represents an upper bound for the probability
that $x,z$ are in one equivalence class and $y,w$ are in another,
while the second summand, $C\,\ta xzyw/^{-\alpha}$,
represents an  upper bound for the probability
that $x,z,y,w$ are all in the same class.
Indeed,  when the
equivalence classes are the
components of the USF,
we will make this heuristic precise in \ref s.usfprop/.

Several examples of random relations that have a stochastic dimension
are described in \ref s.examp/.
The main result of \ref s.usfprop/, \ref t.usfdim/,
asserts that the relation determined by the components
of the USF has stochastic dimension $4$.

\procl d.comp \procname{Composition}  
Let $\LL,\RR\subset\Zd\times\Zd$ be random relations.
The {\bf composition} $\LL\RR$ of $\LL$ and
$\RR$ is the set of all $(x,z)\in\Zd$ such that
there is some $y\in\Zd$ with $x\LL y$ and $y\RR z$.
\endprocl

Composition is clearly an associative operation, that is,
$(\LL \RR)\MM= \LL (\RR\MM)$.  Our main goal in this section
is to prove,

\procl t.comps
Let $\LL,\RR\subset\Zd\times\Zd$ be independent random relations.
Then
$$
\dim_S(\LL\RR)=\min\Bigl\{\dim_S(\LL)+\dim_S(\RR),d\Bigr\},
$$
assuming that $\dim_S(\LL)$ and $\dim_S(\RR)$ exist.
\endprocl


\noindent{\bf Notation.} We write $\phi\lec\psi$
(or equivalently, $\psi \gec \phi$), if
$\phi \le C \psi $ for some constant $C>0$, which
may depend on the laws of the relations considered.
We write $\phi\asymp\psi$ if  $\phi\lec\psi$  and $\phi\gec\psi$.
For $v \in \Zd$ and $0 \le n < N$, define the {\bf dyadic shells}
$$
H_n^N (v):= \{x \in \Zd \st 2^n \le \o{vx} < 2^N\}.
$$

\procl r.restrict
 As the proof will show,
the composition rule of \ref t.comps/ for random relations in
$\Z^d$ is valid
for any graph where the shells
$H_k^{k+1} (v)$
satisfy  $|H_k^{k+1}(v)|\asymp 2^{d k}$.
\endprocl

For sets $V,W \subset\Zd$ let
$$
\dd(V,W):=\min\{\o{vw}:v\in V,\ w\in W\}.
$$
In particular, if $V$ and $W$ have nonempty intersection, then
$\dd(V,W)=1$.
We write $\ta Wx/$ as an abbreviation for $\ta W \cup \{x\} /$.  
\procl l.trees
For every $M>0$, every $x\in\Zd$ and  every $W\subset\Zd$ with
$|W|\le M$, we have
$$
\ta Wx/ \le  \ta W/\,\dd(x,W) \lec \ta Wx/   \,,
$$
where the constant implicit in the $\lec$ notation depends
only on $M$.
\endprocl

\proof
Assume, without loss of generality that $x\notin W$.
The inequality $\ta Wx/ \le \ta W/\,\dd(x,W)$ holds
because given a tree on $W$ we may obtain a tree
on $W\cup\{x\}$ by adding an edge connecting
$x$ to the closest vertex in $W$.
For the second inequality, consider some tree $\tau$
with vertices $W\cup\{x\}$.
Let $W'$ denote the neighbors of $x$ in $\tau$,
and let $u\in W'$ be such that $\o{xu}=\dd(x,W')$.
Let $\tau'$ be the tree on $W$ obtained from $\tau$
by replacing each edge $\{w',x\}$ where $w'\in W' \setminus \{u\}$,
by the edge $\{w',u\}$. (See Figure~\figref{trees}.)  It is easy to verify
that $\tau'$ is a tree.
For each $w'\in W'$, we have
$\o{uw'}\le \o{ux}+\o{xw'}\le 2 \o{xw'}$.
Hence $\ta \tau'/\,\o{ux}\le 2^M \ta \tau/$, and
the inequality $\ta W/\,\dd(x,W)\le 2^M\ta Wx/$ follows.
\qed

\midinsert
\ninepoint
\SetLabels
(.3*.2)$\tau$\\
(.9*.2)$\tau'$\\
\L(.4*.5)$x$\\
\T(.2*.42)$u$\\
\T(.8*.42)$u$\\
\endSetLabels
\centerline{\hskip 0.3in\AffixLabels{\epsfysize=1.8in\epsfbox{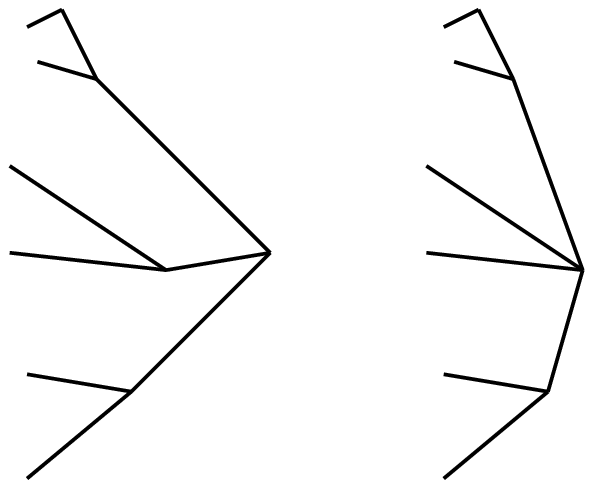}}}
\caption{\figlabel{trees}\enspace
The trees $\tau$ and $\tau'$.}
\endinsert

\procl r.spready
Repeated application of \ref l.trees/ yields that
for any set $\{x_1,\dots,x_n\}$ of $n$ vertices in $\Zd$,
$$
\o{x_1 \dots x_n} \asymp \prod_{i=1}^{n-1}
\dd(x_i, \{x_{i+1},\dots,x_n\}) \,,
\label e.est
$$
where the implied constants depend only on $n$.
\endprocl

Our next goal in the proof of \ref t.comps/ is to establish
\ref e.a/ for the composition $\LL\RR$.  For this, the following
lemma will be essential.

\procl l.comp3
Let $\LL$ and $\RR$ be independent random relations in $\Zd$.
Suppose that $\dim_S(\LL)=d-\alpha$ and $\dim_S(\RR)=d-\beta$ exist, and
denote $\gamma:=\alpha+\beta-d$.
For $u,z \in \Zd$ and $1 \le n \le N$, let
$$
S_{uz}=S_{uz}(n,N):= \sum_{x \in H_n^{N+1} (u)} \I{u\LL x} \I{x \RR z}
\, .
$$
If $\o{uz} < 2^{n-1} $ and $N\ge n$, then
$$
\P[S_{uz}>0] \gec
{\sum_{k=n}^N  2^{-k\gamma} \over \sum_{k=0}^{N}  2^{-k\gamma}} \,.
\label e.2mom
$$
\endprocl

\proof
For $k\ge n$ and $x \in H_k^{k+1}(u)$, we have
$\P[u\LL x] \asymp 2^{-k \alpha}$ (use \ref e.a/ and \ref e.dd/).
Also, for $x \in H_k^{k+1}(u)$,  $k\ge n$,
$$
{1 \over 2} \o{xu}
\le \o{xu} - \o{zu} \le \o{xz} \le \o{xu} + \o{zu} \le 2 \o{xu}
$$
and
$\P[x \RR z] \asymp 2^{-k \beta}$.

Because $\LL$ and $\RR$ are independent and $|H_k^{k+1}(u)|\asymp 2^{d
k}$,
we have
$$
\E[S_{uz}] \asymp \sum_{k=n}^{N} 2^{kd} 2^{-k \alpha} 2^{-k\beta} =
\sum_{k=n}^{N}  2^{-k\gamma} \,. \label e.first
$$
To estimate the second moment, observe that if
$
2\o{uz} \le \o{ux} \le \o{uy}, 
$
then
$$
\o{uy} \le \o{uz}+\o{zy} \le {1 \over 2} \o{uy}+\o{zy} \, ,
$$
 so
$\o{xy}\le\o{xu}+\o{uy} \le 2 \o{uy}  \le 4 \o{yz}$.
Applying the first two inequalities here, \ref e.d/ for $\LL$  and \ref
l.trees/, we obtain that
 $$
\P[u \LL x, u \LL y] \lec
 \oo{ux}{-\alpha} \oo{uy}{-\alpha}+\oo{uxy}{-\alpha} \lec
\oo{ux}{-\alpha}\bigl(
\oo{uy}{-\alpha}+\oo{xy}{-\alpha}\bigr)
\lec
\o{ux}^{\, -\alpha} \o{xy}^{\, -\alpha}.
$$
Similarly, from \ref l.trees/ and the inequality $\o{xy}\le 4 \o{yz}$ above,
it follows that
$$
\P[x \RR z, y \RR z]\lec
\o{xz}^{\,-\beta}\bigl(\o{xy}^{\,-\beta}+\o{yz}^{\,-\beta}\bigr)
\lec \o{xz}^{\, -\beta}  \o{xy}^{\, -\beta} \, .
$$
Since
$$
\E[S_{uz}^2]  \le 2 \sum_{x \in H_n^{N+1} (u)}\
\sum_{y \in H_n^{N+1} (u)}   \II{\o{uy} \ge \o{ux}}  \P[u\LL x, u \LL y]
 \P[x \RR z, y \RR z] ,
$$
we deduce by
breaking the inner sum up into sums over $y \in H_j^{j+1}(x)$
for various $j$ that
$$
\E[S_{uz}^2] \lec \!\!\! \sum_{x \in H_n^{N+1} (u)} \o{ux}^{\, -\alpha}
\o{xz}^{\, -\beta}
 \!\!\! \sum_{y \in H_0^{{N+2}} (x)} \o{xy}^{\, -\alpha-\beta}
\lec \, \E[S_{uz}]
 \sum_{j=0}^{N}  2^{j(d-\alpha-\beta)}
  \,.  \label e.second
$$
Finally, recall that $\gamma=\alpha+\beta-d$, and apply the
Cauchy-Schwarz
inequality
in the form
$$
\P[S_{uz}>0] \ge {\E[S_{uz}]^2 \over \E[S_{uz}^2]} \,.
$$
The estimates \ref e.first/ and \ref e.second/ then yield the assertion
of
the
lemma.
\qed


\procl c.low
Under the assumptions of \ref t.comps/, we have
$\P[u\LL\RR z]\gec \oo{uz}{-\gamma'}$ for all $u,z\in\Zd$,
where $\gamma':=\max\bigl\{0,d-\dim_S(\LL)-\dim_S(\RR)\bigr\}$.
\endprocl

\proof Let $n:=\lceil \log_2 \o{uz}\rceil+2$
and $\gamma:=\alpha+\beta-d$ with $\alpha:=d-\dim_S(\LL)$,
$\beta:=d-\dim_S(\RR)$.
Apply the lemma with $N:=n$ if $\gamma\ne0$ and $N:=2n$ if $\gamma=0$.
\qed

The proof of \ref e.d/ for the composition $\LL\RR$ requires
some further preparation.

\procl l.conv1
Let $\alpha,\beta\in[0,d)$ satisfy $\alpha+\beta> d$.
Let $\gamma:=\alpha+\beta-d$.
Then  for $v,w \in \Zd$,
$$
\sum_{x\in \Zd}
\oo{vx}{-\alpha}\,\oo{xw}{-\beta} \asymp \oo{vw}{-\gamma}  \,.
$$
\endprocl
\proof
Suppose that $2^N \le \o{vw} \le 2^{N+1}$.
By symmetry, it suffices to sum over vertices $x$ such that $\o{xv} \le
\o{xw}$.
For such $x$ in the shell $H_n^{n+1}(v)$, we have by the triangle inequality
$$\oo{xv}{-\alpha}\,\oo{xw}{-\beta} \asymp 2^{-n\alpha }
 2^{- \max\{n,N\}\beta} \,.
$$
Multiplying by $2^{dn}$ (for the volume of the shell) and
summing over all $n$ proves the lemma.
\qed

\procl l.convol
Let $M>0$ be finite and let $V,W\subset \Zd$ satisfy $|V|,|W|\le M$.
Let $\alpha,\beta\in[0,d)$ satisfy $\alpha+\beta > d$.
Denote $\gamma:=\alpha+\beta-d$.
Then
$$
\sum_{x\in\Zd}
\ta Vx/^{-\alpha}\ta Wx/^{-\beta}
\lec
\ta V/^{-\alpha}\,\dd(V,W)^{-\gamma}\,\ta W/^{-\beta}
\le
\ta V W/^{-\gamma}
\,,
\label e.convv
$$
where the constant implicit in the $\lec$ relation may depend
only on $M$.
\endprocl

\proof Using \ref l.trees/ we see that
$$
\ta Vx/^{-\alpha}
\lec
\ta V/^{-\alpha}\, \dd(x,V)^{-\alpha}
\le
\ta V/^{-\alpha}\,\sum_{v\in V}
\oo{vx}{-\alpha}\,
$$
and similarly $\ta Wx/^{-\beta} \lec
 \ta W/^{-\beta} \sum_{w\in W}\,\oo{wx}{-\beta}  $.
Therefore, by \ref l.conv1/,
\begineqalno
\sum_{x\in\Zd}
\ta Vx/^{-\alpha}\ta Wx/^{-\beta}
&\lec
\ta V/^{-\alpha}\, \ta W/^{-\beta} \sum_{v\in V}\sum_{w\in W}
\oo{vw}{-\gamma}
\cr
&\le M^2
\ta V/^{-\alpha}\,\ta W/^{-\beta} \,\dd(V,W)^{-\gamma}\,.
\endeqalno
The rightmost inequality in \ref e.convv/ holds because
$\gamma\le\min\{\alpha,\beta\}$ and given trees on $V$ and
on $W$, a tree on
$V\cup W$ can be obtained by adding an edge $\{v,w\}$ with
$v \in V, w \in W$ and $\o{vw}=\dd(V,W)$ (unless $V\cap W\ne\emptyset$,
in which case $|V\cap W|-1$ edges have to be deleted
to obtain a tree on $V\cup W$).
\qed

\procl l.trip
Let $\alpha,\beta\in[0,d)$ satisfy
$\gamma:=\alpha+\beta-d>0$.  Then
$$
\sum_{a\in\Zd} \oo{ax}{-\alpha}\oo{ay}{-\beta}\oo{az}{-\gamma}
\lec \oo{xyz}{-\gamma}
$$
holds for $x,y,z\in\Zd$.
\endprocl

\proof
Without loss of generality, assume that
$\o{zx} \le \o{zy}$. Consider separately the sum over
$A:=\{a: \o{az} \le {1 \over 2} \o{zx}\}$ and its complement.
For  $a \in A$  we have $\o{ax} \ge   {1 \over 2} \o{zx}$
and  $\o{ay} \ge   {1 \over 2} \o{zy}$.
Therefore,
\begineqalno
\sum_{a \in A} \oo{ax}{-\alpha}\oo{ay}{-\beta}\oo{az}{-\gamma}
&\lec  \oo{zx}{-\alpha} \oo{zy}{-\beta} \sum_{a \in A} \oo{az}{-\gamma}
\cr
& \lec    \oo{zx}{-\alpha}\oo{zy}{-\beta} \oo{zx}{d-\gamma}
\cr &
= \oo{zx}{-\gamma}\oo{zy}{-\gamma}
{\oo{zx}{d-\alpha} \over \oo{zy}{d-\alpha}}
\le \oo{xyz}{-\gamma} \,,
\endeqalno
because of the assumption $\o{zx}\le\o{zy}$ and $\gamma>0$.
Passing to the complement of $A$, we have,
\begineqalno
\sum_{a \in \Zd \setminus A}  \oo{az}{-\gamma} \oo{ax}{-\alpha}
\oo{ay}{-\beta}
&\lec    \oo{zx}{-\gamma}\sum_{a \in \Zd}   \oo{ax}{-\alpha}
\oo{ay}{-\beta}
\cr
&\lec \oo{zx}{-\gamma} \oo{xy}{-\gamma} \le \oo{xyz}{-\gamma} \,,
\endeqalno
where \ref l.conv1/ was used in the next to last inequality.
Combining these two estimates completes the proof of the lemma.
\qed

The following slight extension of \ref l.trip/ will also be
needed.  Under the same assumptions on $\alpha,\beta,\gamma$,
$$
\sum_{a \in \Zd}
\oo{axw}{-\alpha} \oo{ay}{-\beta} \oo{az}{-\gamma}
\lec \oo{xwyz}{-\gamma}.
\label e.ext
$$
This is obtained from \ref l.trip/ by using
$\o{axw}\gec \o{xw}\min\bigl\{\o{ax},\o{aw}\bigr\}$,
which is an application
of \ref l.trees/, and $\o{xw}\o{wyz} \ge \o{xwyz}$, which holds
by the definition of the spread.

\proofof t.comps
If $\dim_S(\LL)+\dim_S(\RR)\ge d$,
then \ref c.low/ shows that
$$
\inf_{x,y\in\Zd}\P[x\LL\RR y]>0\,,
$$
which is equivalent to $\dim_S(\LL\RR)=d$.
Therefore, assume that $\dim_S(\LL)+\dim_S(\RR)< d$.
Let $\alpha:=d-\dim_S(\LL)$,
$\beta:=d-\dim_S(\RR)$ and
$\gamma:=\alpha+\beta-d$.
Since \ref c.low/ verifies \ref e.a/ for the composition $\LL\RR$,
it suffices to prove \ref e.d/ for $\LL\RR$ with $\gamma$ in
place of $\alpha$.
Independence of the relations $\LL$ and $\RR$, together
with \ref e.d/ for $\LL$ and for $\RR$ with $\beta$ in place of $\alpha$
imply
\begineqalno
&
\P[x\LL\RR z,\, y\LL\RR w]
\le \sum_{a,b\in\Zd} \P[x\LL a,\, y\LL b] \P[a\RR z,\, b\RR w]
\cr &\qquad
\lec \sum_{a,b\in\Zd}
\Bigl(\oo{xa}{-\alpha}\oo{yb}{-\alpha}+ \oo{xayb}{-\alpha}\Bigr)
\Bigl(\oo{az}{-\beta}\oo{bw}{-\beta}+ \oo{azbw}{-\beta}\Bigr).
\label e.foursum
\endeqalno
Opening the parentheses gives four sums, which we deal with
separately.
First,
$$
\sum_{a,b\in\Zd}
\oo{xa}{-\alpha}\oo{yb}{-\alpha}
\oo{az}{-\beta}\oo{bw}{-\beta}
\lec \oo{xz}{-\gamma}\oo{yw}{-\gamma},
\label e.ff
$$
by \ref l.conv1/ applied twice.  Second, by two applications
of \ref l.convol/,
\begineqalno
&
\sum_{a\in\Zd}
\sum_{b\in\Zd}
\oo{xayb}{-\alpha} \oo{azbw}{-\beta}
\lec
\sum_{a\in\Zd}
\dd\bigl(\{x,a,y\},\{z,a,w\}\bigr)^{-\gamma}
\oo{xay}{-\alpha}\oo{zaw}{-\beta}
\cr &\qquad\qquad
=
\sum_{a\in\Zd}
\oo{xay}{-\alpha}\oo{zaw}{-\beta}
\lec \oo{xyzw}{-\gamma}.
\label e.ss
\endeqalno
Third, by \ref l.convol/ and \ref e.ext/,
\begineqalno
&
\sum_{a\in\Zd}
\sum_{b\in\Zd}
\oo{xayb}{-\alpha} \oo{az}{-\beta}\oo{bw}{-\beta}
\lec
\sum_{a\in\Zd}
\dd\bigl(w,\{x,y,a\}\bigr)^{-\gamma}\oo{az}{-\beta}
\oo{xay}{-\alpha}
\cr &\qquad\qquad
\le
\sum_{a\in \Zd}
\oo{wa}{-\gamma}\oo{az}{-\beta}
\oo{xay}{-\alpha} +
\dd\bigl(w,\{x,y\}\bigr)^{-\gamma}
\sum_{a\in \Zd}\oo{az}{-\beta} \oo{xay}{-\alpha}
\cr &\qquad\qquad
\lec
\oo{wzxy}{-\gamma}+\dd\bigl(w,\{x,y\}\bigr)^{-\gamma}\oo{zxy}{-\gamma}
\le 2 \oo{wzxy}{-\gamma}.
\label e.tt
\endeqalno
By symmetry, we also have
$$
\sum_{a,b\in\Zd}
\oo{xa}{-\alpha}\oo{yb}{-\alpha} \oo{azbw}{-\beta}
\lec \oo{xywz}{-\gamma}.
$$
This, together with \ref e.foursum/, \ref e.ff/, \ref e.ss/ and
\ref e.tt/ implies  that $\LL\RR$ satisfies the correlation inequality
\ref e.d/  with
$\gamma$ in place of $\alpha$, and completes the proof.
\qed

\bsection {Tail triviality}{s.tail}

Consider the USF on $\Zd$.
Given $n=1,2,\dots$, let $\RR_n$ be the relation consisting of all pairs
$(x,y)\in\Zd\times \Zd$ such that
$y$ may be reached from $x$ by a path which uses no more
than $n-1$ edges outside of the USF.  We show below
that $\RR_1$ has stochastic dimension $4$
(\ref t.usfdim/) and
that $\RR_n$ stochastically dominates the composition of $n$ independent
copies of $\RR_1$ (\ref t.dom/).  By \ref t.comps/,
$\RR_n$ dominates a relation with stochastic dimension $\min\{4n,d\}$.
When $4n\ge d$, this says that $\inf_{x,y\in\Zd}\P[x\RR_n y]>0$.
For \ref t.main1/, the stronger statement that
$\inf_{x,y\in\Zd}\P[x\RR_n y]=1$ is required.
For this purpose, tail triviality needs to be discussed.

\procl d.triv \procname{Tail triviality}
Let $\RR\subset \Zd\times\Zd$ be a random relation
with law $\P$.
For a set $\Lambda \subset\Zd\times\Zd$, let $\field F_\Lambda$ be
the $\sigma$-field generated by the events
$x\RR y$, $(x,y)\in \Lambda$.
Let the {\bf left tail field} $\field F_L(v)$ corresponding to a vertex
$v$
be the intersection of all $\field F_{\{v\} \times K}$
where $K\subset\Zd$ ranges over all subsets
such that $\Zd \smallsetminus K$ is finite.
Let the {\bf right tail field} $\field F_R(v)$ be the intersection
of all $\field F_{K\times \{v\}}$
where $K\subset\Zd$ ranges over all subsets
such that $\Zd \smallsetminus K$ is finite.
Let the {\bf remote tail field}
$\field F_W$ be the intersection of all $\field F_{K_1\times K_2}$,
where $K_1,K_2\subset\Zd$ range over all subsets of $\Zd$
such that   $\Zd \smallsetminus K_1$  and $\Zd \smallsetminus K_2$ are
finite.
The random relation $\RR$ with law $\P$ is said to be {\bf left tail
trivial} if $\P[A]\in \{0,1\}$ for every  $A \in \field F_L(v)$ and
every $v
\in \Zd$.
Analogously, define {\bf right tail triviality} and {\bf remote tail
triviality}.
\endprocl

We will need
the following known lemma, which is a corollary of the main result
of \ref b.weizs/.
Its proof is included for the reader's convenience.
\procl l.tails 
Let $\{\field F_n\}$ and $\{\field G_n\}$ be two
decreasing sequences of complete $\sigma$-fields
in a probability space $(X, \field F, \mu)$,
with $\field G_1$ independent of $\field F_1$,
and let $\field T$ denote the trivial $\sigma$-field,
consisting of events with probability $0$ or $1$.
If $\cap_{n \ge 1} \field F_n =\field T =\cap_{n \ge 1} \field G_n$,
then $\cap_{n \ge 1} (\field F_n \vee \field G_n) = \field T$
as well.
\endprocl

\proof
Let $h\in\bigcap_{n=1}^\infty L^2(\field F_n \vee \field G_n)$.
It suffices to show that $h$ is constant.
Suppose that
$f_n \in L^2(\field F_n)$,  $g_n \in L^\infty(\field G_n)$ and
$f \in L^\infty(\field F_1)$.
Then
$$
\E[f_n g_n f \mid \field G_1]= g_n \E[f_n f \mid \field G_1]
 =g_n \E[f_n f] \qquad a.s.
$$
As the linear span of such products $f_n g_n$ is dense in
$L^2(\field F_n \vee \field G_n)$, it follows that
$$
\E[hf \mid \field G_1] \in L^2(\field G_n) \,.
\label e.obser
$$

By our assumption $\field T=\cap_{n \ge 1} \field G_n$,
we infer from \ref e.obser/ that
$$
 \E[hf \mid \field G_1]=\E[hf]  \,
 \hbox{ for } f \in L^\infty(\field F_1).
\label e.infer
$$
In particular, $\E[h \mid \field G_1]=\E[h]$.
By symmetry,   $\E[h \mid \field F_1]=\E[h]$,
whence
$$
\forall f \in L^\infty(\field F_1),\qquad
 \E[h f]=\E\bigl[f \E[h \mid \field F_1]\bigr]=\E[h] \E[f]\,.
$$
Inserting this into \ref e.infer/,
we conclude that for  $f \in L^\infty(\field F_1)$
and $g \in L^\infty(\field G_1)$ (so that $f$ and $g$ are
independent), we have
$$
\E[h f g]= \E\bigl[g \E[hf \mid \field G_1]\bigr]= \E[h] \E[f]
\E[g] = \E[h] \E[fg] \,.
$$
Thus $h-\E[h]$ is orthogonal to all such products $fg$,
so it must vanish a.s.
\qed

Suppose that $\RR ,\LL\subset\Zd\times\Zd$ are
independent random relations
which have trivial remote tails and trivial left tails,
and have stochastic dimensions.
We do not know whether it follows that
$\LL\RR$ is left tail trivial.
For that reason, we
introduce the notion of the {\bf restricted composition} $\LL\circ\RR$,
which is the relation consisting of all
pairs $(x,z)\in\Zd\times\Zd$ such that
there is some $y\in \Zd$ with
$x\LL y$ and $y\RR z$ and
$$
\o{xz}\le\min\{\o{xy},\o{yz}\}\,.
\label e.restdef
$$

\procl t.restrictedComp
Let $\RR,\LL\subset\Zd\times\Zd$ be independent random relations.
\item{(i)} If $\LL$ has trivial left tail and $\RR$ has
trivial remote tail, then
the restricted composition $\LL\circ \RR$ has trivial left tail.
\item{(ii)} If $\dim_S(\LL)$ and $\dim_S(\RR)$ exist,
then
$$
 \dim_S(\LL\circ\RR)=\min\{\dim_S(\LL)+\dim_S(\RR), \, d\} \,.
$$
\item{(iii)} If $\dim_S(\LL)$ and $\dim_S(\RR)$ exist,
$\LL$ has trivial left tail, $\RR$ has trivial right tail and
$\dim_S(\LL)+\dim_R(\RR)\ge d$, then
$\P[x \LL \circ \RR z]=1$ for all $x,z \in \Zd$.
\endprocl
\proof
\item{(i)} This is a consequence of \ref l.tails/.
\item{(ii)} 
Denote $\gamma':=\max\{0, d-\dim_S(\LL)-\dim_R(\RR)\}$.
The inequality $\P[x \LL\circ\RR z] \gec \oo{xz}{-\gamma'}$
 follows from the proof of \ref c.low/; this concludes the proof
 if $\gamma'=0$. If $\gamma'>0$, then the required upper bound for
 $\P[x \, \LL \circ \RR \,   z,\, y \,\LL \circ \RR \, w]$
 follows from \ref t.comps/,
 since $\LL\circ\RR$ is a subrelation of $\LL\RR$.
\item{(iii)} Let $\field F_n$ (respectively,  $\field G_n$)
be the (completed) $\sigma$-field generated by the events
$u \LL x$ (respectively, $x \RR z$) as $x$ ranges over
$H_n^{\infty}(u)$.
The event
$$
A_n:=\{\exists x\in H_n^{2n}(u)\, : \,  u\LL x,\ x\RR z \}
$$
is clearly in $\field F_n \vee  \field G_n$.
By \ref l.comp3/, there is a constant $C$ such that
$\P[A_n]\geq 1/C>0$, provided that $n$ is sufficiently large.
Let $A=\cap_{k \ge 1} \cup_{n \ge k} A_n$
be the event that there are infinitely many
$x\in\Z^d$ satisfying $u\LL x$ and $x\RR z$.  Then $\P[A]\ge 1/C$.
By \ref l.tails/, $\P[A]=1$.
\qed

\procl c.comp 
Let $m \ge 2$, and let
$\{\RR_i\}_{i=1}^m$ be independent
random relations in $\Zd$ such that $\dim_S(\RR_i)$ exists for each $i
\le
m$.
Suppose that
$$
\sum_{i=1}^{m} \dim_S(\RR_i)  \geq d \, ,
$$
and in addition, $\RR_1$ is left tail trivial,
each of $\RR_2, \dots, \RR_{m-1}$ has a trivial
remote tail, and $\RR_m$ is right tail trivial.
Then  $\P[u\RR_1\RR_2\cdots\RR_m z]=1$
for all $u,z \in \Zd$.
\endprocl
\proof
Let $\LL_1=\RR_1$ and define inductively $\LL_k=\LL_{k-1} \circ \RR_k$
for $k=2,\ldots,m$, so that $\LL_k$
is a subrelation of $\RR_1 \RR_2 \dots \RR_k$.
It follows by induction from \ref t.restrictedComp/ that
 $\LL_k$ has a trivial left tail for each  $k<m$, and
$$
\dim_S(\LL_k) = \min\Bigl\{\sum_{i=1}^{k} \dim_S(\RR_i),d \Bigr\} \,.
$$
By \ref t.restrictedComp/(iii), the restricted composition
$ \LL_m=\LL_{m-1} \circ \RR_m $ satisfies
$$
\P[u\LL_m z]=1 \, \hbox{ for all } u,z \in \Zd \, .
\label e.rest1
$$
\qed

\bsection {Relevant USF properties}{s.usfprop}
Basic to the understanding of the USF is a procedure from
\ref b.BLPS:usf/, that generates
the (wired) USF on any transient graph; it is called
\lq\lq Wilson's method rooted at infinity\rq\rq,
  since it is based on an algorithm from  \ref b.Wilson:alg/
for picking uniformly a spanning tree in a finite graph.
Let $\{v_1,v_2,\dots\}$ be an arbitrary ordering of
the vertices of a transient graph $G$.
Let $X_1$ be simple random walk started
from $v_1$.  Let $F_1$ denote the loop-erasure of $X_1$, which is
obtained by following $X_1$ and erasing the loops as they are
created.  Let $X_2$ be a simple random walk from $v_2$ which
stops if it hits $F_1$, and let $F_2$ be the union of $F_1$
with the loop-erasure of $X_2$.  Inductively, let $X_n$
be a simple random walk from $v_n$, which is stopped
if it hits $F_{n-1}$, and let $F_n$ be the union of
$F_{n-1}$ with the loop-erasure of $X_n$.   Then
$F:=\bigcup_{n=1}^\infty F_n$ has the distribution of the (wired) USF on
$G$.
The edges in $F$  inherit the
orientation from the loop-erased walks creating them, and hence
$F$ may be thought of as an oriented forest.
Its distribution does not depend on the ordering
chosen for the vertices.
See \ref b.BLPS:usf/ for details.

We say that a random set $\A$ {\bf stochastically dominates} a
random set $\Q$ if there is a coupling $\mu$ of $\A$
and $\Q$ such that $\mu[\A \supset \Q]=1$.
\procl t.dom \procname{Domination}
Let $F,F_0,F_1,F_2,\dots,F_m$ be independent samples of the
wired USF in the graph $G$.
Let $\C(x,F)$ denote the vertex set of the component of $x$ in $F$.
 Fix a distinguished vertex $\v$ in $G$, and write
$\C_0 =\C(\v,F)$. For $j \ge 1$, define inductively
$\C_{j}$ to be the union of all vertex components of $F$
that are contained in, or adjacent to, $\C_{j-1}$.
Let  $\Q_0=\C(\v,F_0)$.  For $j \ge 1$, define inductively
$\Q_{j}$ to be the union of all vertex components of $F_j$
that intersect $\Q_{j-1}$.
Then $\C_m$ stochastically dominates $\Q_m$.
\endprocl

\proof
For each $R>0$ let
$B_R:=\bigl\{v\in\Zd\st |v|<R\bigr\}$, where $|v|$
is the graph distance from $\v$ to $v$.
Fix $R>0$.
Let $B_R^W$ be the graph obtained from $G$ by collapsing
the complement of $B_R$ to a single vertex, denoted $\vi$.
Let $F',F'_0,\dots,F'_m$ be independent samples of the
uniform spanning tree (UST) in $B_R^W$.
Define $F^*$ to be $F'$ without the edges incident to $\vi$,
and define $F^*_i$ for $i=0,\ldots,m$ analogously.

 Write
$\C_0^* =\C(\v,F^*)$. For $j \ge 1$, define inductively
$\C_{j}^*$ to be the union of all vertex components of $F^*$
that are contained in, or adjacent to, $\C_{j-1}^*$.

Let  $\Q_0^*=\C(\v,F_0^*)$.  For $j \ge 1$, define inductively
$\Q_{j}^*$ to be the union of all vertex components of $F_j^*$
that intersect $\Q_{j-1}^*$.


We show by induction on $j$ that $\C^*_j$ stochastically dominates
$\Q^*_j$. The induction base where $j=0$ is obvious.
For the inductive step, assume that $0 \le j<m$ and there is a
coupling $\mu_j$ of $F'$ and $(F'_0,F'_1,\dots,F'_j)$ such that
$\C^*_j\supset \Q^*_j$
holds $\mu_j$-a.s.
Let $H$ be a set of vertices in
$B_R$ such that $\C_j^*=H$ with
positive probability. Let $\widehat{H^c}$ denote the subgraph
of $B_R^W$ spanned by the vertices
$B_R\cup\{\vi\} \sm H$.
Let $S_H$ be $F'\cap \widehat{H^c}$ conditioned
on $\C^*_j=H$.   Since $\C_j^*$
is a union of components of $F^*$, it is clear that $S_H$ has the
same distribution as a UST on $ \widehat{H^c}$.  Therefore,
by the negative association property of uniform spanning trees (see the discussion
following Remark~5.7 in \ref b.BLPS:usf/),
$S_H$ stochastically dominates $F'_{j+1}\cap \widehat{H^c}$
(conditional on $\C^*_j=H$).
Hence, we may extend $\mu_j$ to a coupling
$\mu_{j+1}$ of $F'$ and $(F'_0,F'_1,\dots,F'_{j+1})$ such that
$\C^*_j\supset \Q^*_j$ and  $H=\C_j^*$
satisfies $F^* \cap \widehat{H^c} \supset F^*_{j+1}\cap \widehat{H^c}$
a.s.\ with respect to $\mu_{j+1}$. In such a coupling,
we also have $\C^*_{j+1}\supset \Q^*_{j+1}$ a.s.

This completes the induction step, and proves that $\C_m^*$
stochastically dominates $\Q^*_m$.
Observe that $\C^*_m$ converges weakly to $\C_m$
and $\Q^*_m$
converges weakly to $\Q_m$, as $R\to\infty$.
This follows from the fact that
$F$ stochastically dominates $F^*$ (see \ref b.BLPS:usf/,
Cor. 4.3(b)) and
$F^* \to F$ weakly as $R\to\infty$.
The theorem follows.
\qed

\procl t.usfdim \procname{Stochastic dimension of USF}
Let $F$ be the USF in $\Zd$, where $d\geq 5$.  Let
$\UU$ be the relation
$$
\UU:=\bigl\{(x,y)\in\Zd\times\Zd\st
\hbox{$x$ and $y$ are in the same component of $F$}\bigr\}
\,.
$$
Then $\UU$ has stochastic dimension $4$.
\endprocl
\proof  We use Wilson's method rooted at infinity
and a technique from \ref b.LPS/.
For $x,z \in \Z^d$, consider independent simple random walk paths
$\{X(i)\}_{i \ge 0}$ and $\{Z(j)\}_{j \ge 0}$
that start at $x$ and $z$ respectively.
By running Wilson's method rooted at infinity starting with $x$ and
then starting another random walk from $z$, we see
that $\P[z\UU x]$ is equal to the probability that the walk $Z$
intersects the loop-erasure of the walk $X$.
Given $m\in\N$, denote by $\{L_m(i)\}_{i=0}^{q_m}$
the path obtained from loop-erasing $\{X(i)\}_{i =0}^m$.
Define
\begineqalno
 \tau_L(m)&:=\inf\Bigl\{k \in\{0,1,\dots,q_m\}
\st L_m(k) \in \{X(i)\}_{i \ge m}\Bigr\}\, ,
\cr
 \tau_L(m,n)&:=\inf\Bigl\{k \in\{0,1,\dots,q_m\}
\st L_m(k) \in \{Z(j)\}_{j \ge n}
\Bigr\}\,.
\cr
\endeqalno
Consider the indicator variables
$
I_{m,n}:=\II{X(m)=Z(n)}
$
and
$$
J_{m,n}:=I_{m,n} \II{\tau_L(m,n) \le \tau_L(m)}\,.
$$
Observe that on the event $J_{m,n}=1$, the path
$\{Z(j)\}_{j \ge 0}$
intersects the loop-erasure of  $\{X(i)\}_{i =0}^\infty$
at $L_m\bigl(\tau_L(m,n)\bigr)$.
Hence $\P[x\UU z]= \P[\sum_{m,n} J_{m,n}>0]$.
Given $X(m)=Z(n)$, the law of
$\bigl\{X(m+j)\bigr\}_{j\ge 0}$
is the same as the law of
$\bigl\{Z(n+j)\bigr\}_{j\ge 0}$.
Therefore,
$$
\P[\tau_L(m,n) \le \tau_L(m) \mid X(m)=Z(n)] \ge 1/2 \,,
$$
so that $\E[I_{m,n}]\ge \E[J_{m,n}] \ge \E[I_{m,n}]/2$.
Let
\begineqalno
&\Phi=\Phi_{xz}:=\sum_{m=0}^\infty\   \sum_{n=0}^\infty I_{m,n} \,;
\cr
&\Psi=\Psi_{xz}:=\sum_{m=0}^\infty\  \sum_{n=0}^\infty J_{m,n} \,.
\cr
\endeqalno
Then
\begineqalno
\E[\Psi] & \asymp   \sum_{m=0}^\infty\   \sum_{n=0}^\infty \E[I_{m,n}]
\cr
&=
 \sum_{m=0}^\infty\   \sum_{n=0}^\infty\   \sum_{v \in \Zd} \P[X(m)=v]
\,\P[Z(n)=v]
\cr
&=  \sum_{v \in \Zd} G(x,v) G(z,v) \, ,
\cr
\endeqalno
where $G$ is the Green function for simple random walk in $\Zd$.
Since $G(x,v) \asymp \o{xv}^{\, 2-d}$ (see, e.g., \ref b.Spitzer/),
we infer that (see \ref l.conv1/)
$$
\E[\Psi] \asymp \sum_{v \in \Zd}  \o{xv}^{\, 2-d} \o{zv}^{\, 2-d} \asymp
\o{xz}^{\, 4-d} \,.
\label e.spitz
$$
The second moment calculation for $\Phi$ is classical.
Again, the relation
$G(x,v) \asymp \o{xv}^{\, 2-d}$ is used:
\begineqalno
 \E[\Psi^2] & \le \E[\Phi^2]
\cr
&
\le  \sum_{y,w \in \Zd} G(x,y)G(y,w)[G(z,y)G(y,w)+G(z,w)G(w,y)]
\cr
& \lec \sum_{y \in \Zd} G(x,y)G(z,y)+ \sum_{y\in \Zd}
G(x,y)\sum_{w \in \Zd}G(y,w)^2 G(z,w)
\cr
& \lec \sum_{y \in \Zd} \o{xy}^{\, 2-d}\,\o{zy}^{\,2-d}+ \sum_{y\in \Zd}
\o{xy}^{\, 2-d}\sum_{w \in \Zd}\o{yw}^{\, 4-2d}\,\o{zw}^{\,2-d}
\,.
\endeqalno
The proof of~\ref l.conv1/ gives
$$
\sum_{w \in \Zd}\o{yw}^{\, 4-2d}\,\o{zw}^{\,2-d}\lec \o{yz}^{\,2-d}\,.
$$
Hence
$$
 \E[\Psi^2]
 \lec
\o{xz}^{\, 4-d}+
\sum_{y\in \Zd} \o{xy}^{\,2-d}\,\o{yz}^{\,2-d} \lec \o{xz}^{\, 4-d} \, .
$$
Therefore,
$$
\P[x\UU z]\ge
\P[\Psi_{xz}>0] \ge {\E[\Psi_{xz}]^2 \over  \E[\Psi_{xz}^2]} \gec
\o{xz}^{\,
4-d}  \,.
$$
This verifies \ref e.a/ for $\UU$
with $d-\alpha=4$.

To verify \ref e.d/, we must bound
$\P[x \UU z, y \UU w]$.
Let $X,Z,Y,W$ be simple random walks starting from $x,z,y,w$,
respectively.
Let $\hat X$ denote the set of vertices visited by $X$, and similarly
for $Z,Y$ and $W$.
Then we may generate the USF by first using the random
walks $X,Z,Y,W$.
On the event $x\UU z \wedge y\UU w\wedge\neg x\UU y$,
we must have $\hat X\cap \hat Z\ne \emptyset$
and $\hat Y\cap\hat W \ne \emptyset$.
Hence,
\begineqalno
\P[ x\UU z \wedge y\UU w\wedge\neg x\UU y]
&
\le
\P[\hat X\cap\hat Z\ne\emptyset]\, \P[\hat Y\cap\hat W\ne\emptyset]
\cr &
\le
\sum_{a\in\Zd} G(x,a)G(z,a)\sum_{b\in\Zd}G(y,b)G(w,b)
\cr &
\lec \o{xz}^{\,4-d} \o{yw}^{\,4-d}.
\label e.nottogether
\endeqalno
On the other hand, the event
$ x\UU z \wedge y\UU w\wedge x\UU y$ is the event
$U(x,y,z,w)$ that $x,y,z,w$ are all in the same USF tree.

By \ref t.spread/ below,
$\P[U(x,y,z,w))]\lec \oo{xyzw}{4-d}$.
This together with \ref e.nottogether/ gives
$$
\P [ x\UU z \wedge y\UU w]\le
\P[ x\UU z \wedge y\UU w\wedge\neg x\UU y]+\P[U(x,y,z,w)]
\lec \oo{xz}{4-d} \oo{yw}{4-d}+\oo{xyzw}{4-d}\,,
$$
which verifies \ref e.d/ with $\alpha=d-4$,
and completes the proof.
\qed

\procl t.spread
For any finite set $W \subset \Zd$, denote by $U(W)$ the event
that all vertices in $W$ are in the same USF component. Then
$$
\P[U(W)] \lec \oo{W}{4-d} \,,
\label e.together
$$
where the implied constant depends only on $d$ and the cardinality of
$W$.
\endprocl

\proof
When $|W|=2$, say $W=\{x,z\}$, \ref e.together/ follows from \ref
e.spitz/,
since $\P[ x\UU z] = \P[\Psi>0]\le \E[\Psi]$.
We proceed by induction on $|W|$.
For the inductive step, suppose that
$|W| \ge 3$. For  $x,y \in W$ denote by $U(W;x,y)$ the intersection of 
 $U(W)$ with the event that the path connecting $x,y$ in the USF
is edge-disjoint from the oriented USF paths connecting the vertices in
$V:=W\setminus\{x,y\}$ to infinity.
 By considering all the possibilities for the vertex $z \in \Zd$
where the oriented USF paths from $x$ and $y$ to $\infty$ meet,
and running Wilson's method rooted at infinity starting
with random walks from the vertices in $V\cup\{z\}$ and following
by random walks from $x$ and $y$
we obtain
$$ 
\P[U(W;x,y)]  \lec \sum_{z \in \Zd}
\P\bigl[U(Vz)\bigr]\oo{zx}{2-d}
\oo{zy}{2-d} \,,
\label e.findz
$$
where $Vz:=V\cup\{z\}$.
 By the induction hypothesis and \ref l.trees/,
$$
 \P[U(Vz)] \lec  \oo{Vz}{4-d} \lec
 \oo{V}{4-d} \sum_{v \in V}  \oo{vz}{4-d} \,.
$$
Inserting this into \ref e.findz/ and applying \ref l.trip/
with $\alpha=\beta=d-2$, yields
\begineqalno
\P[U(W;x,y)] & \lec \oo{V}{4-d} \sum_{v \in V}
  \sum_{z \in \Zd} \oo{vz}{4-d} \oo{zx}{2-d} \oo{zy}{2-d}
\cr
& \lec  \oo{V}{4-d} \sum_{v \in V} \oo{vxy}{4-d}  \lec \oo{W}{4-d}
\,,
\endeqalno
where the last inequality follows from the fact that
 for every $v \in V$, the union of a tree on $\{v,x,y\}$ and a tree
on $V$ is a tree on $W=V\cup\{x,y\}$.
Finally, observe that
$$ 
U(W)  = \bigcup_{(x,y)} U(W;x,y) \,,
$$
where $(x,y)$ runs over all pairs of distinct vertices in $W$.
Indeed, on $U(W)$, denote by $m(W)$
 the vertex where all oriented USF paths based in $W$
meet, and pick $x,y$ as the pair of vertices in
$W$ such that their oriented USF paths meet farthest
(in the intrinsic metric of the tree)
from $m(W)$.
Consequently,
$$
\P[U(W)]  \le \sum_{(x,y)} \P[U(W;x,y)] \lec \oo{W}{4-d} \,.
$$
\qed

\procl r.sharp
The estimate in \ref t.spread/ is tight, up to constants, i.e.,
for any finite set $W \subset \Zd$,
$$
\P[U(W)] \gec \oo{W}{4-d} \,,
\label e.together2
$$
where the implicit constant may depend only on $d$ and the cardinality
of $W$.
Since we will not need this lower bound, we omit the proof.
\endprocl

\procl t.usftriv \procname{Tail triviality of the USF relation}
The relation $\UU$ of \ref t.usfdim/ has trivial left, right
and remote tails.
\endprocl

In \ref b.BLPS:usf/ it was proved that the tail of the (wired or free)
USF on every infinite graph is trivial.  However, this is not the same
as
the tail triviality of the relation $\UU$.
Indeed, if the underlying graph $G$ is a regular tree of degree greater
than $2$, then the relation $\UU$ determined by the (wired)
USF in $G$ does not have a trivial left tail.

\proof
The theorem clearly holds when $d\le 4$, for then $\UU=\Zd\times\Zd$
a.s.  Therefore, restrict to the case $d>4$. Let $F$ be the USF
in $\Zd$. We start with the remote tail.
Fix $\xx\in\Zd$, $r>0$, and let $\mmS=\{K_1 \subset F, K_2 \cap F =\emptyset\}$
be a cylinder event where $K_1, K_2$ are sets of edges in
a ball $B(\xx,r)$ of radius $r$ and center $x$.
For $R>0$, let $\Y_{R}$ be the union of all one-sided-infinite
simple paths in $F$ that start at some vertex $v\in\Zd\sm B(x,R)$.
By \ref b.BLPS:usf/, Theorem 10.1, a.s.\ each component of $F$ has one
end.
(This implies that $\Y_R$ is a.s.\ the same as the union of all
oriented USF paths starting outside of $B(x,R)$.)
Therefore, there exists a function $\phi(R)$ with $\phi(R)\to\infty$
as $R\to\infty$ such that
$$
\lim_{R\to\infty}\P[\Y_R\cap B(\xx,\phi(R))\ne\emptyset]=0\,.
\label e.noint
$$
Let $\tilde\Y_R:=\Y_R$ if $\Y_R\cap B(\xx,\phi(R))=\emptyset$,
and $\tilde\Y_R:=\emptyset$ otherwise.
Note that $\tilde\Y_R$ is a.s.\ determined
by
$F\setminus B(\xx,\phi(R))$, or more precisely by the set of
edges in $F$ with both endpoints outside $B(x,\phi(R))$.
By tail triviality of $F$ itself, it therefore follows that
$$
\lim_{R\to\infty}
\EBig{\bigl|\P[\mmS|\tilde\Y_R]-\P[\mmS]\bigr|}=0
\,.
\label e.purim
$$
For $R>r$, on the event
$\Y_{R}\cap B\bigl(\xx,\phi(R)\bigr)=\emptyset$, we have
$\P[\mmS \mid \tilde \Y_R]=\P[\mmS\mid\Y_R]$.
Hence, by~\ref e.noint/,
$$
\lim_{R\to\infty} \EBig{\bigl|\P[\mmS \mid \Y_R]-\P[\mmS \mid \tilde
\Y_R]\bigr|}=0
\,.
$$
With~\ref e.purim/ this gives
$$
\lim_{R\to\infty}
 \EBig{\bigl|\P[\mmS \mid \Y_R]-\P[\mmS]\bigr|}=0
\,.
\label e.yrs
$$
Since the remote tail of $\UU$ is determined by $\Y_{R}$ for every $R$,
the remote tail is independent of $\mmS$, whence it is trivial.

Now consider the left tail of $\UU$.
Let $A$ be an event in $\field F_L(\xx)$.
Let $\mmS=\{K_1 \subset F, K_2 \cap F =\emptyset\}$  as above,
where $K_1, K_2$ are sets of edges in $B(\xx,r)$.
To establish the triviality of $\field F_L(x)$,
it is enough to show that $A$ and $\mmS$ are independent.

Let $\Y'_R$ be the union of all the one-sided-infinite simple
paths in $F$ that start at some vertex in the outer boundary of
$B(x,R)$. By Wilson's method rooted at infinity we may construct
$F$ by first choosing $\Y'_R$, then the path in $F$, $\gar$ say,
from $x$ to $\Y'_R$ and after that the paths starting at points
outside $\Y'_R \cup \gar$. Let $\zeta_R$ be the endpoint of
$\gar$ on $\Y'_R$. Recall that $\gar$ is obtained by
loop erasure of a simple random walk path starting from $x$ and
stopped when it hits $\Y'_R$. Given $\Y'_R$,
 the paths in $F$ to $\Y'_R$ from the vertices
in the complement of $\Y'_R \cup B(x,R)$ are conditionally independent of
$\gar$. Therefore the conditional
distribution of $\zeta_R$ given $\Y_R$ is just the harmonic measure on
$\Y'_R$ for a simple random walk started at $x$. Given $\Y_R$,
we shall write $\mu_y(z; \Y_R)$ for the probability that a
simple random walk started at $y$ first hits $\Y_R$ at $z$.
Similarly, for a set $W$ of vertices we write
$\mu_y(W;\Y_R):=\sum_{z\in W} \mu_y(z;\Y_R)$.

It is clear that the pair $(\Y_R, \zeta_R)$ determines for which
$z \notin B(x,R)$ the relation $x \UU z$ holds. Therefore, the indicator
function of $A$ is measurable with respect to the pair
$(\Y_R, \zeta_R)$.
Given $\Y_R$, let $\AR$ denote the set of vertices $z\in\Y_R$ such that
$A$ holds if $\zeta_R=z$.  Then
$$
\P[A] = \Pbig{\zeta_R\in \AR} = \Ebig{\mu_x(\AR; \Y_R)} \,.
\label e.hk1
$$

A very similar relation holds for $\P[\mmS, A]$. Instead of
choosing $\gar$ immediately after $\Y'_R$, we can first
determine whether $\mmS$ occurs after choosing $\Y'_R$ and then
choose $\gar$. Again when $\Y'_R$ is given, the paths from
points outside $\Y'_R$ and outside $B(x,R)$ are not influenced
by $\mmS$, so that $\P[\mmS \mid \Y'_R] = \P[\mmS \mid \Y_R]$. We further
remind the reader of the following fact (see, e.g., \ref b.BLPS:usf/).
Suppose that $G=(V,E)$ is a finite graph and
$E_1,E_2\subset E$.
Let $T'$ be the (set of edges of the) UST in $G$.
Then $T'$, conditioned on $T'\cap (E_1\cup E_2)= E_1$
(assuming this event has positive probability),
is the union of $E_1$ with the set
of edges of a UST on the graph obtained
from $G$ by contracting the edges in $E_1$ and deleting
the edges in $E_2$. 
Let $H$ be the graph obtained from $\Zd$ by contracting
the edges in $K_1$ and deleting the edges in $K_2$.
By first choosing $\Y'_R$ and continuing with 
Wilson's algorithm we see that
the conditional law of $F \setminus \Y_R$ given $\Y_R$, is the 
law of a UST on the finite graph obtained from $\Zd$ 
by gluing all vertices in $\Y_R$ to a single vertex. If we further
condition on the occurrence of $\mmS$, then we should also contract
the edges in $K_1$ and delete the edges in $K_2$.
Therefore, conditionally on $\Y_R$ and
the occurrence of $\mmS$, the path $\gar$ has the distribution of the
loop-erasure of a simple random walk on $H$, starting at $x$ and
stopped when it hits $\Y'_R$. If we write $\mu^H_y(z; \Y_R)$ for
the probability that a
simple random walk on $H$ started at $y$ first hits $\Y_R$ at
$z$, 
then we obtain analogously to \ref e.hk1/ that
$$
\P[\mmS,A] = 
\P[\mmS,\zeta_R\in\AR] =
 \EBig{\P[\mmS \mid \Y_R]\, \mu^H_x(\AR; \Y_R)}
\,.
$$
In view of \ref e.yrs/ this gives
$$
\Bigl|\P[\mmS,A] - \P[\mmS]\,  \Ebig{\mu^H_x(\AR; \Y_R) }\Bigr|
\le \EBig{\bigl|\P[\mmS \mid \Y_R]-\P[\mmS]\bigr|}
\mathop{\longrightarrow}\limits_{R\to\infty}0\,.
\label e.hk3
$$

Given $\epsilon>0$, let $R_1>r$ be such that a simple random walk
started from any vertex outside $B(x,R_1)$ has probability
at most $\epsilon$ to visit $B(x,r)$.
Since every bounded harmonic function in $\Zd$ is constant, there exists
$R_2>R_1$ such that any harmonic function $u$ on $B(\xx,R_2)$
that takes values in $[0,1]$, satisfies the Harnack inequality
$$
\sup_{y,z \in B(\xx,R_1+1)} \bigl|u(y)-u(z)\bigr|<\epsilon \,.
\label e.harnack
$$
(See, e.g., Theorem 1.7.1(a) in \ref b.Lawler/.)
In particular, when $\Y_R \cap B(x,R_2) = \emptyset$, we can apply this to 
the harmonic function $y\mapsto \mu_y(\AR; \Y_R)$ to obtain
$$
\sup_{y, y' \in B(x, R_1+1)} \bigl|\mu_y(\AR;\Y_R) - \mu_{y'}(\AR;\Y_R)\bigr| < \epsilon\,
\quad\hbox{provided that}\quad
\Y_R \cap B(x,R_2) = \emptyset\,
 .
\label e.hk4
$$
We need to show that $\mu_\xx(\AR,\Y_R)$ is close to $\mu^H_\xx(\AR,\Y_R)$.
To this end, note that by the definition of $R_1$,
when a simple random walk on $H$ from $\xx$ first exits $B(\xx,R_1)$,
 it has probability at most $\epsilon$ to
revisit $B(\xx,r)$.  On the event that it does not,
it may be considered as a random walk in $\Z^d$ which does
not visit $B(\xx,r)$.
By \ref e.hk4/ it follows that
$$
\bigl|\mu_\xx^H(\AR;\Y_R)- \mu_\xx(\AR;\Y_R)\bigr|\le 2\epsilon\,,
\label e.top5
$$
on the event $\Y_R \cap B(x,R_2) = \emptyset$.
Finally take $R_3>R_2$ sufficiently large that
$$
\forall R\ge R_3,\qquad\P[\Y_{R}\cap B(\xx,R_2)\ne\emptyset]<\epsilon
\,.
$$
{}From \ref e.top5/ it follows that for $R \ge R_3$
$$
\Bigl|\Ebig{\mu^H_x(\AR; \Y_R)} -\Ebig{ \mu_x(\AR; \Y_R)} \Bigr| < 3 \epsilon\,.
$$
Combined with \ref e.hk3/ and \ref e.hk1/ this shows that for
sufficiently large $R$
$$
\bigl|\P[\mmS, A] - \P[\mmS]\,\P[A]\bigr| < 4\epsilon.
$$
As $\epsilon>0$ was arbitrary,
we conclude that $A$ and $\mmS$ are independent,
which proves that $\UU$ has trivial left tail.  By symmetry, the right
tail is trivial too.
\qed

\procl r.scope
The above proof of left-tail triviality for $\UU$
is valid for the wired USF in any transient graph $G$ such that
there are no non-constant bounded harmonic functions and
a.s.\ each component of the USF has one end.
Only the one end property is needed for triviality of the remote tail.
(For recurrent graphs, the USF is a.s.\ a tree, whence obviously
the USF relation has trivial left, right and remote tails.)
\endprocl

The following lemma will be needed in the proof of the last statement
of \ref t.main1/.

\procl l.sb
Let $D\subset\Zd$ be a finite connected set with a connected
complement, and denote by $\widehat{D^c}$ the
subgraph of $\Zd$ spanned by the vertices
in $\Zd\setminus D$.
Let $F$ be the USF on $\Zd$, and 
denote by $\Gamma_D$ the event that there
 are no oriented edges in $F$ from $D^c$ to $D$.
  Then the distribution of $F\cap \widehat{D^c}$ conditioned on
$\Gamma_D$, is the same as the distribution of the wired USF in the 
graph $\widehat{D^c}$.
\endprocl

\proof
We first consider a finite version of this statement.
Let $G=(V,E)$ be a finite connected graph, $\rho \in V$ a
distinguished vertex, and $D \subset V\setminus\{\rho\}$.
Let $G_-$ be the subgraph of $G$ spanned by $V \setminus D$.
Denote by $S(G)$ the set of spanning trees in $G$,
and let $\Gamma_D^*$ be the set of 
$t \in S(G)$ such that for every $w \in D^c$, the path
 in $t$ from $w$ to $\rho$ is disjoint from $D$.

\noindent{\bf Claim:} {\it Let $T$ be a uniform spanning tree (UST) in $G$.
Assume that $\Gamma_D^*\ne\emptyset$.
Then  conditioned on
$T \in \Gamma_D^*$, the edge set
$T \cap G_-$ is a UST in $G_-$. } 

To prove the claim, fix two trees $t_1, t_2 \in S(G_-)$,
and for every $t \in \Gamma_D^*$ that contains $t_1$,
define $t':=(t\setminus t_1) \cup t_2$. For each vertex $v \in V$
there is a path  in $t'$ from $v$ to $\rho$, that uses edges
of $(t\setminus t_1)$ until it reaches $D^c$, and then uses edges of 
$t_2$.
Note that $t$ cannot
contain any edge of $t_2 \setminus t_1$, because $t_1$ plus any
edge of $\widehat D^c$ outside $t_1$ contains a circuit. Thus
$t$ and $t'$ have the same number of edges. It follows that $t'
\in S(G)$ and moreover, $t' \in \Gamma_D^*$. The map $t
\mapsto t'$, is a bijection, because
$t' \mapsto (t' \setminus t_2) \cup t_1$ is its inverse.
This shows that $t_1$ and $t_2$ have
the same number of extensions to spanning trees of $G$ that are
in $\Gamma_D^*$. Moreover any $t \in \Gamma_D^*$ is an
extension of some $t_1 \in S(G_-)$. In other words, we have
established the claim.

The lemma follows by considering the uniform spanning
forest on $\Zd$ as a weak limit of uniform spanning trees in finite
subgraphs (with wired complements),
where $\rho$ is chosen as the wired vertex.
\qed

\procl r.fincomp
If a finite connected set
$D\subset\Zd$ has a connected complement,
then the lemma above implies
that a.s., every component of
the wired USF in $\Zd \sm D$ has one end.
The proof of \ref t.usftriv/ can then be adapted to show that the
wired USF relation in $\Zd \sm D$ has trivial remote, right, and left, tail
$\sigma$-fields. Note that this can also be inferred from \ref r.scope/.
 To verify that every bounded harmonic function on $\Zd \sm D$ is constant,
we recall that the existence of  nonconstant bounded harmonic functions on
a graph is equivalent to the existence of two disjoint sets of
vertices $A,B$ such that with positive probability the random walk
eventually stays in $A$, and the same holds for $B$.
(If such  $A$ and $B$ exist, then define the harmonic function
$h(v)$ as the probability that simple random walk started from $v$ eventually
stays in $A$. If $h$ is a bounded harmonic function with $\sup h=1$,
$\inf h=0$, say, then we may take $A:=h^{-1}([0,1/4])$ and
$B:=h^{-1}([3/4,1])$.)
This criterion is clearly uneffected by the removal of $D$ from $\Zd$.
\endprocl

\bsection{Proofs of main results}{s.mainproofs}
\proofof t.main1
Denote $m=\lfloor {{d-1}\over 4} \rfloor$.
By applying \ref c.comp/ to $m+1$ independent copies
of the USF relation, and invoking Theorems \briefref t.usfdim/
and \briefref t.usftriv/, we infer that a.s., every two vertices in
$\Zd$ are related by the composition of these $m+1$ copies.
Now \ref t.dom/ yields the upper bound
$\max\bigl\{N(x,y)\st x,y\in\Zd\bigr\} \le m$ a.s.

For the final assertion of the theorem,
it suffices to show that for every $x,y\in\Zd$ and every $r>0$,
the event $\Xi(x,y,r,m)$ that there
is a path in $\Zd \setminus B_r$ from $T(x)$ to $T(y)$
with at most $m$ edges outside $F$, satisfies 
$$
\P[\Xi(x,y,r,m)]=1 \,.
\label e.goalie
$$

Let $\Delta_d(r)$ denote the
collection of finite, connected sets $D \subset \Zd$
such that $B_r \subset D$ and $D^c$ is connected.
 For $D \in \Delta_d(r)$,  denote by $\Gamma_D$ the event that
there are no oriented edges in $F$ from $D^c$ to $D$.
We claim that for every $r>0$,
$$
\P\Bigl[\bigcup_{D \in \Delta_d(r)} \Gamma_D \Bigr]=1 \,.
\label e.ID
$$
Indeed, let $\aB_r$ denote the set of vertices $v \in \Zd$ such that
the oriented path from $v$ to infinity in $F$ enters $B_r$.
Since a.s.\ each of the USF components has one end, 
$\aB_r$ and its complement are connected, and $\aB_r$ is a.s.\ finite. 
Thus $\Gamma_{\aB_r}$ occurs and \ref e.ID/ holds.

Therefore, to prove \ref e.goalie/, it suffices to show
that 
$$
\forall D \in \Delta_d(r), \qquad
\P\bigl[\Xi(x,y,r,m) \mid \Gamma_D\bigr]=1 \,.
\label e.goalie2
$$
Fix $D \in \Delta_d(r)$,
let $G_-$ be the subgraph of $\Zd$ spanned by the vertices
in $\Zd\setminus D$, and let $F_-$ denote the wired
USF on $G_-$.
By \ref l.sb/, conditioned on $\Gamma_D$,
the law of $F\cap G_-$ is the same as that of $F_-$.
Consequently, it suffices to show that a.s.\ every pair of
vertices in $G_-$ is connected by a path in $G_-$ with at
most $m$ edges outside of $F_-$.
The USF relation on $G_-$ has trivial left, right
and remote tails, by \ref r.fincomp/.  It is also
clear that this relation has stochastic dimension $4$.
Therefore, the above proof for $\Zd$ applies also to
$G_-$, and shows that every pair of
vertices in $G_-$ is connected by a path in $G_-$ with at
most $m$ edges outside $F_-$.
Since $r$ is arbitrary and $D$ is
an arbitrary set in $ \Delta_d(r)$,
this proves the last assertion of the theorem,
and also completes the proof for $5 \le d \le 8$.

For the case $d>8$, it remains only to prove the lower bound on
 $\max N(x,y)$.  This is done in the following proposition.
\qed

\procl p.lower
The inequality
$$
\P\bigl[N(x,y)\le k\bigr]\lec\oo{xy}{4k+4-d}
$$
holds for all $x,y\in\Zd$ and all $k\in\N$.
\endprocl

\proof
The proposition clearly holds for $k\ge\bigl \lfloor (d-1)/4\bigr
\rfloor$.
So assume that $k< \bigl\lfloor (d-1)/4 \bigr\rfloor$.
Consider a sequence $\{x_j,y_j\}_{j=0}^k$ in $\Zd$ with
$x_0=x$, $y_k=y$ and $|y_j-x_{j+1}|=1$ for $j=0,1,\dots,k-1$.
Let $A=A(x_0,y_0,\dots,x_k,y_k)$ be the event
that $y_j\in T(x_j)$ for $j=0,1,\dots,k$ and
$y_j\notin T(x_i)$ if $j\ne i$.
Observe that
$$
\bigl\{N(x,y)=k\bigr\}\subset\bigcup A(x_0,y_0,\dots,x_k,y_k)\,,
$$
where the union is over all such sequences $x_0,y_0,\dots,x_k,y_k$.
To estimate the probability of $A$,
we run Wilson's algorithm rooted at infinity, begining with
random walks started at $x_0,y_0,\dots,x_k,y_k$.
For $A$ to hold, for each $j=0,\dots,k$, the random
walk started at $x_j$ must intersect the path of the random
walk started at $y_j$.
Hence, as in \ref e.spitz/,
$$
\P[A]\lec \prod_{j=0}^k \oo{x_jy_j}{4-d}\,.
\label e.abnd
$$
Since $x_{j+1}$ and  $y_j$ are adjacent, this implies
$$
\P\bigl[N(x,y)=k\bigr]
\lec \sum \prod_{j=0}^k \oo{x_jx_{j+1}}{4-d}\,,
\label e.N
$$
where the sum extends over all sequences
$\{x_j\}_{j=0}^{k+1}\subset\Zd$ such that
$x_0=x$ and $x_{k+1}=y$.

We use \ref l.conv1/ repeatedly in \ref e.N/,
first to sum over $x_{k}$, then over
$x_{k-1}$, and so on, until $x_1$, and obtain
$\P\bigl[N(x,y)=k\bigr]\lec \oo{xy}{4k+4-d}$, which completes the proof.
\qed

%

\proofof t.main2
The case $d\le 4$ is from \ref b.Pemantle/.
The case $d=5,6,7,8$ follows from \ref t.main1/.
For $d \ge 9$ and fixed $x,y \in \Zd$ and $r>0$
we infer from \ref e.abnd/ that
\begineqalno
&
\Pbig{\exists v \in T(x), \,\exists w \in T(y) : |v| \ge n, \,  |w| \ge
n, \,
|v - w| \le r, \, T(x) \ne T(y)}
\cr
&\qquad
\le \sum_{v,w} \Pbig{A(x, v, y,w)}
\lec \sum_{v,w} \oo{xv}{4-d}\oo{yw}{4-d}
 \, ,
\endeqalno
where the sums are over $v,w$ satisfying
$|v|, |w| \ge n$ and $|v -w| \le r$.  The latter sum is bounded by
$C(x,y,r)\, n^{8-d}$, by \ref l.conv1/.
\qed

\bsection{Further Examples and Remarks}{s.examp}

We now present some examples of relations having a stochastic
dimension, and several questions.

\procl x.triv
Fix some $\alpha\in (0,d)$.
Let $\Rel R\subset\Zd\times\Zd$ be the random relation such that
$\{x\Rel R y\}_{\{x,y\}\subset \Zd}$ are independent events,
and $\P[x\Rel R y]= \oo{xy}{-\alpha}$.
Then $\Rel R$ has stochastic dimension $d-\alpha$. This relation
corresponds
to {\it long range percolation} in $\Zd$, where the degree of every
vertex
is infinite and the infinite cluster contains all vertices of $\Zd$.
\ref t.comps/ implies that the intrinsic diameter of this graph
is $\lceil {d \over d-\alpha} \rceil$ almost surely.
\endprocl

\procl x.walk
Let $d\ge 3$, and
for every $x\in\Zd$, let $S^x$ be a simple random walk
starting from $x$, with $\{S^x\}_{x\in\Zd}$ independent.
Let $x\Rel S y$ be the relation $\exists n\in\N\ \ S^x(n)=y$.
It is easily seen that $\Rel S$ has stochastic dimension $2$.
{}From Theorems \briefref t.usfdim/, \briefref c.comp/
and \briefref t.usftriv/,
we infer that simple random walk a.s.\ intersects
each component of an independent USF iff $d \le 6$.
\endprocl

\procl x.shake
For every $x\in\Zd$, let $S^x$ be as above.
Define a relation $\Rel Q$ where $x\Rel Qy$ iff $\exists n\in\N\ \ S^x(n)=S^y(n)$.
Then $\Rel Q$ has stochastic dimension $2$ (provided $d\ge 2$).
The proof is left to
the reader.
\endprocl

\procl x.coal
Fix $d\ge 2$.
For every $x\in\Zd$, let $S^x_*(\cdot)$ be a continuous time simple
random
walk starting from $x$,  where particles
walk independently until they meet; when two particles meet, they
coalesce,
and the resulting particle carries both labels.
Let $x\Rel W y$ be the relation $\exists t>0\ \ S^x_*(t)=y$.
Then $\Rel W$ has  stochastic dimension $2$.
\endprocl



\procl x.othergreen \procname{USF generated by other walks}
Let $\{ S_n \}_{n \ge 0}$ be a symmetric random walk on $\Z^d$,
i.e., the increments $S_n - S_{n-1}$ are symmetric i.i.d.\ $\Zd$-valued
variables.
Say that this random walk has {\bf Greenian index} $\alpha$
if the Green function
$
G_S(x,y) := \sum_{n=0}^{\infty} \P [S_n =y \vert S_0 =x ]
$
satisfies
$
G_S(x,y)
\asymp
\oo{xy}{\alpha-d}
$
 for all $x,y \in \Z^d$.
It is well known that random walks with bounded increments of mean zero
have Greenian index 2, provided that $d \ge 3$;
for $d=3$, the boundedness assumption may be relaxed
 to finite variance (see \ref b.Spitzer/).
\ref b.Williamson/ shows that  for $0 < \alpha < \min\{2,d\}$, if
a random walk $\{S_n\}$ in $\Z^d$ with mean zero increments satisfies
$\P [S_n-S_{n-1} = x] \sim c| x | ^ {-d-\alpha}$
(i.e., the ratio of the two expressions here tends to
1 as $|x| \to \infty$), then this random walk
 has Greenian index $\alpha$.
We can use Wilson's method rooted at infinity to generate
the wired spanning forest (WSF)
based on a symmetric random walk with Greenian index $\alpha$.
(Note, this will not be a subgraph of the usual graph
 structure on $\Z^d$.)
The arguments of the preceding section extend to show that the
relation determined by the components of this WSF
will have stochastic dimension $2\alpha$.
\endprocl


\procl x.min \procname{Minimal Spanning Forest}
 Let the edges of $\Zd$ have i.i.d.~weights
$w_e$, which are uniform random variables in $[0,1]$.
The {\bf minimal spanning forest} is  the subgraph of
$\Zd$ obtained by removing every edge that has the
maximal weight in some cycle;  see, e.g., \ref B.Newman-Stein/.
\endprocl

\procl g.minimal
Let $x\Rel M y$ if $x$ and $y$ are in the same
minimal spanning forest component.
Then $\Rel M$ has stochastic dimension $8$ in $\Zd$ if
$d\ge 8$.  This is a variation on a conjecture of \ref B.Newman-Stein/.
\endprocl

\procl r.dims
It is natural to ask how the notion of stochastic dimension is related
to other notions of dimension for subsets of a lattice. In this
direction, we state the following.
Let $\LL\subset\Zd\times\Zd$ be a random relation.
Denote  $\Gamma_o=\Gamma_o(\LL):=\{z : o\LL z\}$; 
if $\LL$ is an equivalence relation,
then $\Gamma_o(\LL)$ is the equivalence class of the origin $o$,
but we do not assume this.
\endprocl

\procl p.dims
Suppose that $\LL \subset\Zd\times\Zd$ has
stochastic dimension $\alpha \in (0,d) $. Denote by $\eta_n$ the number of vertices
in $\Gamma_o(\LL) \cap B(o,n)$. Then 
\item{(i)} The laws of $\eta_n/n^\alpha$ are tight, but do not converge weakly
to a point mass at $0$. 
\item{(ii)} With positive probability,
$\limsup_n {\log \eta_n \over \log n}=\alpha$.

\item{(iii)} If, moreover, $\LL$ has a trivial left tail, then
$$\limsup_n {\log \eta_n \over \log n} =\alpha \, \hbox{  a.s. }
$$
\endprocl
\proof
(i) The definition of stochastic dimension,
and a standard calculation, 
yield that $\E \eta_n \asymp n^\alpha$ and 
$\E \eta_n^2 \asymp n^{2\alpha}$. These estimates imply 
the asserted tightness, and also that 
$$
\inf_n \P(\eta_n \ge \E\eta_n/2) >0
\label e.Kahane
$$
(see, e.g., Kahane (1985), p.\ 8).
\item{(ii)} The bound $\E \eta_n \lec n^\alpha$ implies that
$\sum_k \eta_{2^k}/(2^{k\alpha} k^2 )< \infty$ a.s. Thus, monotonicity of
$\eta_n$ yields that $\limsup_n {\log \eta_n \over \log n} \le \alpha$ a.s.
The lower bound follows from \ref e.Kahane/.
\item{(iii)} This follows from (ii) and the definition of the left tail.
\qed

When $\LL$ is the USF relation and $d>4$, we have $\alpha=4$, 
and the limsup in (iii)
can be replaced by a limit (we omit the proof).
With more work, it can be shown that
the ``discrete Hausdorff dimension'' 
(in the sense of Barlow and Taylor (1992))
of the USF component $\Gamma_o$ is also 4, almost surely.

A harder task is to prove existence of the {\it scaling limit } for the USF
in dimension $d>4$, and show that this limiting object has
Hausdorff dimension 4 almost surely.

An easily stated problem in this direction, is
to show that 
$$\P[\Gamma_o \cap B(x,k) \neq \emptyset] \asymp (k/|x|)^{d-4}
 \label e.scale
$$
for $k<|x|$. At present, we can only establish such an estimate up
to a power of $\log x$. Note that \ref e.scale/, if true,
depends on properties of the USF beyond its
stochastic dimension; the analog of \ref e.scale/ fails for the relation
considered in \ref x.triv/.

\procl x.quaint 
As noted in the introduction, for $5 \le d \le 8$, the USF defines a
translation invariant
random partition of $\Zd$ with infinitely many connected components,
 where any two components come within distance $1$ infinitely often.
A partition with these properties exists for any $d \ge 3$:
In the hyperplane $H:=\{x \in \Zd : x_d=0\}$, let $L_0$
be the lattice consisting of all vertices with all coordinates even,
and let $L$ be a random translate of $L_0$ in $H$, chosen uniformly
among the $2^{d-1}$ possibilities.
Let $F_0$ be the union of all lines perpendicular to $H$
which contain a vertex of $L$, and complete $F_0$ to a spanning forest
$F$
using Wilson's algorithm.
\endprocl

\procl x.noperc
Consider the random relation $\PP\subset\Z^2\times\Z^2$, where $x\PP y$
iff $x$ and $y$ are in the same connected component of
Bernoulli bond percolation at the critical parameter $p=p_c=1/2$.
It is known that a.s.\ all the components will be finite.  If
$\LL,\RR\subset\Z^2\times\Z^2$ and
$\bigl|\{y:x\LL y\}\bigr|<\infty$ and $\bigl|\{y:x\RR y\}\bigr|<\infty$
for all $x\in\Z^2$, then
 $\{y:x\LL\RR y\}$ is finite too.  By inductively applying  \ref
t.comps/
to independent copies of $\PP$, it follows that $\PP$ cannot have
a positive stochastic dimension.  However,  the Russo-Seymour-Welsh
Theorem implies that $\P[x\PP y]\ge \o{xy}^{-\beta}$
for some $\beta<\infty$.
\endprocl

\procl r.other
Our results on the USF in $\Zd$ extend readily to other
Cayley graphs of polynomial growth.  Let $G$ be such a graph,
and let $o$ be a vertex of $G$.
By \ref b.Gromov/, $G$ satisfies $|B(o,r)|\asymp r^d$ for some integer
$d$ (see also \ref b.WoessRWIGG/, Lemma 3.14 and Theorem 3.16). 
Consequently, we may find an $a>2$ such that
$\bigl|B(o,a^{n+1})\setminus B(o,a^n)\bigr|\asymp a^{nd}$ for $n\in\N$.
Thus, in our above proofs, we may work with shells based on powers of
$a$, in place of the shells $H_n^N$.
By \ref b.HeSalCo/, the Green function estimates which we used in $\Zd$
apply to $G$.   Moreover, any bounded harmonic function on $G$ is a constant
(see, e.g., \ref b.KaiVer/).
This implies, in particular, that the free and the wired
USF on $G$ are the same, by Thm.~7.3 of \ref b.BLPS:usf/ (so it is
clear what we mean by the USF on $G$).
By Theorem 10.1 of \ref b.BLPS:usf/ a.s.\ each component of the USF on
$G$ has one end, unless $G$ is quasi-isometric to $\Z$ (in which case
$G$ is recurrent, $d=1$, the USF is a tree a.s.\ and our results clearly hold).
With these basic facts established, the proofs go through just as for $\Zd$,
with the case $d\le 4$ handled by Cor.~9.6 of \ref b.BLPS:usf/.
\endprocl

\noindent{\bf Acknowledgment} We thank Russ Lyons for
useful discussions, and Jim Pitman for a reference.
H.~Kesten thanks the Hebrew University for its hospitality and
for appointing him as Schonbrunn Visiting Professor during the
Spring of 1997, when this project was started.
We also thank the referee for helpful comments.

\citeinc{Lyons:usfsurvey}
\citeinc{Wilson:alg}

\def\endreferences{\bigskip\bigskip\endgroup}
\def\startreferences{
 \vskip0pt plus.3\vsize \penalty -150 \vskip0pt
 plus-.3\vsize \bigskip\bigskip \vskip \parskip
 \begingroup\baselineskip=12pt\frenchspacing
 \bibliographytitle
 \vskip12pt\parindent=0pt
 \def\and{{\rm and}}
 \def\em{\it}
 \def\newblock{\hskip .11em plus.33em minus.07em}
 \def\bibauthor##1{{\sc ##1}}
 \def\bibitem[##1]##2
 {\htmlanchor{##2}{}\RefLabel{##2}[##1]\hangindent=.8cm\hangafter=1}
 }
\def\endreferences{\bigskip\bigskip\endgroup}

\ifundefined{bibstylemodification}\relax\else\bibstylemodification\fi
\startreferences

\bibitem[Barlow and Taylor (1992)]
\bibauthor{Barlow, M. \and{} Taylor, S.J.} (1992).
\newblock  Defining fractal subsets of $\Z^d$.  
\newblock {\em Proc. Lon. Math. Soc.} {\bf 64}, 125--152.

\bibitem[BLPS (2001)]{BLPS:usf}
\bibauthor{Benjamini, I., Lyons, R., Peres, Y., \and{} Schramm, O.}
(2001).
\newblock Uniform spanning forests.
\newblock {\em Ann. Prob.} {\bf 29}, 1--65.

\bibitem[Gromov (1981)]{Gromov}
\bibauthor{Gromov, M.} (1981).
\newblock Groups of polynomial growth and expanding maps
\newblock {\em Inst. Hautes \'Etudes Sci. Publ. Math.} {\bf 53}, {53--73}.

\bibitem[Hebisch and Saloff-Coste (1993)]{HeSalCo}
\bibauthor{Hebisch, W. \and{} Saloff-Coste, L.} (1993).
\newblock Gaussian estimates for Markov chains and random walks on
              groups.
\newblock {\em Ann. Prob.} {\bf 21}, 673--709.

\bibitem[Kahane (1985)]
\bibauthor{Kahane, J.-P.} (1985). 
\newblock {\em Some Random Series of Functions.} 
\newblock Second edition,  Cambridge University Press, Cambridge.

\bibitem[Kaimanovich and Vershik (1983)]{KaiVer}
\bibauthor{Ka\u{\i} manovich, V.A. \and{} Vershik, A.M.},
\newblock Random walks on discrete groups: boundary and entropy.
\newblock {\em Ann. Probab.} {\bf 11}, {457--490}.

\bibitem[Lawler (1991)]{Lawler}
\bibauthor{Lawler, G.F.} (1991).
\newblock {\em Intersections of Random Walks}.
\newblock Birkh\"auser Boston Inc., Boston, MA.

\bibitem[Lyons (1998)]{Lyons:usfsurvey}
\bibauthor{Lyons, R.} (1998).
\newblock A bird's-eye view of uniform spanning trees and forests.
\newblock In {\em Microsurveys in discrete probability (Princeton, NJ,
1997)},
  pages 135--162. Amer. Math. Soc., Providence, RI.

\bibitem[LPS (2002)]{LPS}
\bibauthor{Lyons, R., Peres, Y., \and{} Schramm, O.} (2001).
\newblock Markov chain intersections and the loop-erased walk.
\newblock {\em Ann. Inst. H. Poincar\'e, Probab. et Stat.}, to appear.

\bibitem[Newman and Stein (1994)]{Newman-Stein}
\bibauthor{Newman, C.M. \and{} Stein, D.} (1994).
\newblock Spin-glass model with dimension-dependent ground state.
\newblock {\em Phys. Rev. Lett.} {\bf 72}, 2286--2289.

\bibitem[Pemantle (1991)]{Pemantle}
\bibauthor{Pemantle, R.} (1991).
\newblock Choosing a spanning tree for the integer lattice uniformly.
\newblock {\em Ann. Probab.} {\bf 19}, 1559--1574.

\bibitem[Spitzer (1976)]{Spitzer}
\bibauthor{Spitzer, F.} (1976).
\newblock {\em Principles of Random Walk}.
\newblock Springer-Verlag, New York, second edition.
\newblock Graduate Texts in Mathematics, Vol. 34.

\bibitem[von Weizs{\"a}cker (1983)]{weizs}
\bibauthor{von Weizs{\"a}cker, H.} (1983).
\newblock Exchanging the order of taking suprema and countable
intersections of $\sigma$-algebras.
\newblock {\em Ann. Inst. H. Poincar\'e, Probab. et Stat.} {\bf 19},
91--100.

\bibitem[Williamson (1968)]{Williamson}
\bibauthor{Williamson, J. A.} (1968).
\newblock Random walks and Riesz kernels.
\newblock{\em Pacific J. Math.} {\bf 25}, 393--415.

\bibitem[Wilson (1996)]{Wilson:alg}
\bibauthor{Wilson, D.B.} (1996).
\newblock Generating random spanning trees more quickly than the cover
time.
\newblock In {\em Proceedings of the Twenty-eighth Annual ACM Symposium
  on the Theory of Computing (Philadelphia, PA, 1996)}, pages
296--303, New York.  ACM.

\bibitem[Woess (2000)]{WoessRWIGG}
\bibauthor{Woess, W.} (2000).
\newblock {\em Random walks on infinite graphs and groups}, volume 138 of {\em
  Cambridge Tracts in Mathematics}.
\newblock Cambridge University Press, Cambridge.

\endreferences

\filbreak
\begingroup
\eightpoint\sc
\parindent=0pt\baselineskip=10pt

\emailwww{itai@wisdom.weizmann.ac.il}
{http://www.wisdom.weizmann.ac.il/\string~itai/}
\medskip
\emailwww{kesten@math.cornell.edu}{}
\medskip
\emailwww{peres@stat.berkeley.edu}
{http://www.stat.berkeley.edu/\string~peres/}
\medskip
\emailwww{schramm@microsoft.com}
{http://research.microsoft.com/\string~schramm/}
\endgroup
\bye